\documentclass[11pt]{article}

\usepackage{amsmath,amssymb,amsthm,setspace,graphicx,array,tikz,url,blkarray,wasysym,easyReview}
\usepackage[text={6.5in,8.4in}, top=1.3in]{geometry}
\newtheorem{theorem}{Theorem}[section]

\newtheorem{lemma}[theorem]{Lemma}
\newtheorem{corollary}[theorem]{Corollary}

\newtheorem{claim}[theorem]{Claim}

\begin{document}
	
\title{
		Acyclic subgraphs of digraphs with high chromatic number
}
	
\author{
	Raphael Yuster
	\thanks{Department of Mathematics, University of Haifa, Haifa 3103301, Israel. Email: raphy@math.haifa.ac.il\;.}
}
	
\date{}
	
\maketitle
	
\setcounter{page}{1}
	
\begin{abstract}
	For a digraph $G$, let $f(G)$ be the maximum chromatic number of an acyclic subgraph of $G$.
	For an $n$-vertex digraph $G$ it is proved that $f(G) \ge n^{5/9-o(1)}s^{-14/9}$ where $s$ is the bipartite independence number of $G$, i.e., the largest $s$ for which there are two disjoint $s$-sets of vertices with no edge between them.
	This generalizes a result of Fox, Kwan and Sudakov,
	who proved this for the case $s=0$ (i.e., tournaments and semicomplete digraphs).
	Consequently, if $s=n^{o(1)}$, then $f(G) \ge n^{5/9-o(1)}$
	which polynomially improves the folklore bound $f(G) \ge n^{1/2-o(1)}$.
	As a corollary, with high probability, all orientations of the random $n$-vertex graph with edge probability $p=n^{-o(1)}$ (in particular, constant $p$, hence almost all $n$-vertex graphs)
	satisfy $f(G) \ge n^{5/9-o(1)}$. Our proof uses a theorem of Gallai and Milgram that together with several additional ideas, essentially reduces to the proof of Fox, Kwan and Sudakov.
	
\vspace*{3mm}

\end{abstract}

\section{Introduction}

Throughout this paper, the chromatic number $\chi(G)$, the independence number $\alpha(G)$, and the bipartite independence number $\alpha^*(G)$\,\footnote{Recall that $\alpha(G)$ is the largest $k$ for which there is an empty $k$-vertex set and that $\alpha^*(G)$ is the largest $k$ for which there are two disjoint $k$-vertex sets with no edges between them.} of a digraph $G$ are inherited from the corresponding parameters of its underlying undirected graph.
An {\em orientation} (a.k.a. oriented graph) is a digraph without directed cycles of length $2$. A {\em tournament} is an orientation of a complete graph. A digraph is {\em acyclic} if it has no directed cycles.

Graph colorings and orientations are closely related topics whose origins date back to the
classical theorem of Gallai–Hasse–Roy–Vitaver \cite{gallai-1968,hasse-1965,roy-1967,vitaver-1962}
asserting that every orientation of a $k$-chromatic graph has a directed path with $k$ vertices.
Many variations and related results connecting these concepts can be found in the survey \cite{havet-2013}.
A very general variant of the Gallai–Hasse–Roy–Vitaver Theorem was conjectured by Burr \cite{burr-1980},
stating that any $(2k-2)$-chromatic digraph contains every oriented tree on $k$ vertices.
The special case of Burr's conjecture for tournaments is known as
Sumner's conjecture which has been resolved, for $k$ sufficiently large, by K\"uhn, Mycroft, and Osthus
\cite{KMO-2011}. The best bound on the chromatic number of digraphs which guarantees the containment of every oriented tree on $k$ vertices
is $O(k^{1.5})$ which was recently obtained by Bessy, Gon{\c{c}}alves, and Reinald \cite{BGR-2024}.

Burr's conjecture is known to hold for acyclic digraphs. In fact, it was proved by
Addario-Berry, Havet, Sales, Reed, and Thomass\'e \cite{addario+-2013} that every $k$-chromatic acyclic digraph 
contains every oriented tree on $k$ vertices. It would therefore be of interest to prove that
digraphs have acyclic subgraphs that retain the original chromatic number to some nontrivial extent.
Indeed, this problem was formulated in \cite{addario+-2013}.
For a digraph $G$, let $f(G)$ denote the maximum chromatic number of an acyclic subgraph of $G$
and let $f(k)$ be the minimum of $f(G)$ taken over all $k$-chromatic digraphs.
It is a folklore argument (see below) that $f(G) \ge \sqrt{\chi(G)}$ (whence $f(k) \ge \sqrt{k}$), and it was proved in \cite{addario+-2013}
that $f(k) \ge 1+\lfloor \sqrt{k-1} \rfloor$. Trivially, $f(3)=2$ and it is known that $f(4)=3$ \cite{NY-2019,shapira-2022}. The special case of this parameter for the case of tournaments was considered
by Nassar and the author in \cite{NY-2019}. Let $g(k)$ be the maximum integer such that every $k$-vertex tournament has an acyclic subgraph with chromatic number $g(k)$. It was conjectured in \cite{NY-2019}
that $g(k) = \omega(\sqrt{k})$, namely that the folklore lower bound can be improved by more than a constant factor. This conjecture has recently been resolved by Fox, Kwan, and Sudakov \cite{FKS-2021} who proved by an
ingenious argument that $g(k)=\Omega(k^{5/9-o(1)})$. In that same paper, they have shown that the exponent $5/9$
cannot be improved beyond $3/4$. Returning back to $f(k)$, this result suggests that $f(k)=\omega(\sqrt{k})$ as well.
Our main result shows that this is the case in a dense regime.
 
\begin{theorem}\label{t:main}
	Let $G$ be an $n$-vertex digraph with $\alpha^*(G) \le s$.
	Then, $f(G) \ge n^{5/9-o(1)}s^{-14/9}$. In particular, if $\alpha^*(G)=n^{o(1)}$, then
	$f(G) \ge n^{5/9-o(1)}$.
\end{theorem}
\noindent
Theorem \ref{t:main} extends the result of Fox, Kwan, and Sudakov which is the case $s=0$.
Notice that if $\alpha^*(G) \le s$, then $\alpha(G) \le 2s+1$ and thus $\chi(G) \ge n/(2s+1)$.
So, Theorem \ref{t:main} significantly improves upon the folklore $\sqrt{\chi(G)}$ lower bound
for graphs with relatively small bipartite independence number.
Theorem \ref{t:main} is particularly interesting in the setting of {\em random graphs}.
Let ${\cal G}(n,p)$ denote the  Erd\H{o}s–R\'enyi random $n$-vertex graph with edge probability $p$.
It is well known that if $p = n^{-o(1)}$, then $G \sim {\cal G}(n,p)$ almost surely\footnote{As usual, ``almost surely'' means with probability tending to $1$ as $n$ tends to infinity.} has
$\chi(G) =n^{1-o(1)}$ and $\alpha^*(G)=n^{o(1)}$ (see \cite{bollobas-1988} for sharp estimates when $p$ is constant). We therefore have:
\begin{corollary}\label{coro:main}
	Let $p = n^{-o(1)}$ and $G \sim {\cal G}(n,p)$. Then almost surely,
	every orientation $\vec{G}$ of $G$ has $f(\vec{G}) \ge n^{5/9-o(1)}$. Using $p=1/2$ we have that almost all $n$-vertex graphs $G$ have the property that every orientation $\vec{G}$ of $G$ has $f(\vec{G}) \ge n^{5/9-o(1)}$.
\end{corollary}

To prove Theorem \ref{t:main}, we first use a theorem of Gallai and Milgram on the path-covering number of digraphs with given independence number. This is then used, together with some additional ideas, to effectively
reduce the proof to the same framework of the proof of Fox, Kwan, and Sudakov. We do, however, need to (sometimes extensively) modify their arguments in locations where their proof relies on the fact that every pair of vertices is adjacent, which is not the case in our setting. In Section \ref{sec:lem} we prove several lemmas needed for our proof.
In Section \ref{sec:fks} we revisit the proof of Fox, Kwan and Sudakov, state variants of their
Lemmas 2.2 (which is Lemma \ref{l:l22} here) and 2.3 (which is Lemma \ref{l:l23} here) and prove
Lemma \ref{l:l23}. Lemma \ref{l:l22}, which is the main technical lemma, is proved in Section \ref{sec:l22}. Throughout the paper, we ignore rounding issues as long as they do not affect the asymptotic claims.

We end the introduction recalling the aforementioned folklore lower bound. Clearly, in every red/blue edge coloring
of a graph $G$, $\chi(G_r) \cdot \chi(G_b) \ge \chi(G)$ where $G_r$ and $G_b$ respectively denote the red and blue subgraphs. For every digraph $G$ and every linear order of its vertices, we can color red the edges going from left to right and color blue the edges going from right to left, observing that the obtained $G_r$ and $G_b$ are acyclic. Hence $\max(\chi(G_r),\chi(G_b)) \ge \sqrt{\chi(G)}$.
Nevertheless, despite the similarity between acyclic orientations and red/blue edge colorings, note that an analogous result such as Theorem \ref{t:main} is, in general, false for red/blue edge colorings. To see this, consider the complete $n/s$-partite
graph with each part of order $s$; clearly in this case $\alpha^*(G) \le \alpha(G)=s$.
Further partition the parts into $\sqrt{n/s}$ buckets, each bucket consisting of$\sqrt{n/s}$ parts. Now color blue all the edges with endpoints in parts belonging to the same bucket
and color red all remaining edges. It is immediate to verify that $G_r$ and $G_b$ are both $\sqrt{n/s}$ chromatic.

\section{A few helpful lemmas}\label{sec:lem}

We shall require the following result of Gallai and Milgram \cite{GM-1960} (see also \cite{diestel-2005}).
\begin{lemma}\label{l:gm}{\rm\cite{GM-1960}}
	Let $G$ be a digraph with $\alpha(G)=s$. Then there is a set of $s$ pairwise vertex-disjoint directed paths of $G$ that cover all its vertices.
\end{lemma}

For a permutation $\pi$ of the vertices of a digraph $G$, let $G_\pi$ denote the spanning subgraph of $G$
consisting of all edges $(u,v)$ for which $\pi(u) < \pi(v)$. Notice that $G_\pi$ is acyclic and that  $G_{\pi^{rev}}$ consists of the edges going {\em against} $\pi$.

For the rest of this paper, all digraphs are assumed to be orientations. Clearly this assumption does not change the chromatic number, independence number, or bipartite independence number.
For a subset $S$ of vertices of an oriented graph $G$, let $G[S]$ denote the induced subgraph of $G$ on $S$.
We say that $S$ is an {\em acyclic $k$-set} if $|S|=k$ and $G[S]$ is acyclic.
\begin{lemma}\label{l:prob}
	Let $G$ be a $k$-vertex orientation with $\alpha(G) \le s$. For a random permutation $\pi$, the probability that $G_\pi$ has no edges is at most $(se/k)^k$.
\end{lemma}
\begin{proof}
	By Lemma \ref{l:gm}, $G$ has a set of $s$ pairwise vertex-disjoint paths $P_1,\ldots,P_s$ that cover each vertex.
	Let $p_i$ denote the number of vertices of $P_i$ (possibly $p_i=0$) and observe that $\sum_{i=1}^s p_i = k$.
	Notice that for $G_\pi$ to have no edges, all vertices of $P_i$ must
	appear in the reverse order they appear in $P_i$, which, for a random $\pi$, occurs with probability
	$1/p_i!$. As the $P_i$'s are pairwise vertex-disjoint, we have that the probability that
	$G_\pi$ has no edges is therefore at most the reciprocal of $\prod_{i=1}^s p_i!$.
	Finally note that by logarithmic convexity of $x^x$,
	$$
	\prod_{i=1}^s p_i! \ge \prod_{i=1}^s \left( \frac{p_i}{e}\right)^{p_i} \ge \left( \frac{\sum_{i=1}^s p_i}{se}\right)^{\sum_{i=1}^s p_i}=\left( \frac{k}{se}\right)^{k}\;.
	$$
\end{proof}
\begin{corollary}\label{coro:prob}
	If $G$ is an $n$-vertex orientation with $\alpha(G) \le s$ and $G$ has fewer than $(k/se)^k$ acyclic $k$-sets, then there is a permutation $\pi$ such that $\alpha(G_\pi) \le k$. In particular, $\chi(G_\pi) \ge n/k$.
\end{corollary}
\begin{proof}
	Consider a random $\pi$. For a set $S$ of $k$ vertices, the
	probability that $S$ is an independent set in $G_\pi$  is either $0$ if $S$ is not acyclic, or, by Lemma \ref{l:prob}, at
	most $(se/k)^k$ if $S$ is acyclic. Hence, if $G$ has fewer than $(k/se)^k$ acyclic $k$-sets, then by the union bound, $\alpha(G_\pi) \le k$ with positive probability.
\end{proof}

We need an additional simple lemma which guarantees that in any orientation with $\alpha(G) \le s$, there is always a vertex whose in-degree is not too small and whose out-degree is not too small.

\begin{lemma}\label{l:deg}
	Let $G$ be an $m$-vertex orientation with $\alpha(G) \le s$. Then $G$ has a vertex whose in-degree and out-degree are both at least $m/4s-1/2$.
\end{lemma}
\begin{proof}
	Let $t=m/4s-1/2$. Assume the claim is false and partition $V(G)$ into parts $V_1$ and $V_2$ where all vertices of $V_1$ have out-degree smaller than $t$ and all vertices of $V_2$ have in-degree smaller than $t$. Without loss of generality, $|V_1| \ge m/2$. Consider the underlying undirected graph of $G[V_1]$. Since $\alpha(G[V_1]) \le s$, there is a subgraph of this undirected graph
	with minimum degree at least $|V_1|/s-1 \ge m/2s-1$. Returning to the directed version, there must be some vertex in this subgraph whose out-degree is not smaller than its in-degree, hence its out-degree is at least $(m/2s-1)/2=t$, contradicting the definition of $V_1$.
\end{proof}

\section{Generalization of the Fox-Kwan-Sudakov Theorem}\label{sec:fks}

Analogous to the notion of $q$-almost transitive in \cite{FKS-2021}, we say that an orientation is {\em $q$-almost-acyclic} if there is a permutation $\pi$ such that for every vertex $v$, its undirected degree in $G_{\pi^{rev}}$ is at most $q$.
In other words, at most $q$ edges incident to $v$ are oriented against $\pi$.
Using Lemma \ref{l:prob}, Corollary \ref{coro:prob}, Lemma \ref{l:deg}, and the definition of $q$-almost-acyclic, we can now generalize the result from \cite{FKS-2021} to apply to orientations with a given bipartite independence number.

We use the same flow of their proof. In some places, their proof generalizes verbatim, but in other places, their proof relies on the fact that there is an edge between every pair of vertices (namely, tournaments) while in our case, this is not so, hence we need to adjust.

The following two lemmas are analogous to Lemmas 2.2 and 2.3 in \cite{FKS-2021}, respectively.

\begin{lemma}\label{l:l22}
	Let $k = n^{\Omega(1)}$ and $k \ge s^7$. If an $n$-vertex orientation $G$ with $\alpha(G) < s$ and $\alpha^*(G) < s$ has at least $(k/se)^k$ acyclic $k$-sets, then it has a set of $n' \ge k^{2-o(1)}/s$ vertices inducing a $q$-almost-acyclic subgraph $G'$
	where $q =  n's^6/k^{1-o(1)}$.
\end{lemma}

\begin{lemma}\label{l:l23}
Suppose $n^{1/3}s^{-2/3} \le q \le n^{1-\Omega(1)}$. An $n$-vertex $q$-almost-acyclic orientation $G$ with
$\alpha(G) < s$ has a permutation $\pi$ with $\alpha(G_\pi) \le n^{1/4+o(1)}q^{1/4}s^{3/2}$, so $\chi(G_\pi) \ge n^{3/4-o(1)}q^{-1/4}s^{-3/2}$.
\end{lemma}

\begin{proof}[Proof of Theorem \ref{t:main}]
	Since for any graph $G$ it holds that $\alpha(G) \le 2\alpha^*(G)+1$, and since the theorem's outcome stays intact up to $n^{o(1)}$ (in particular, if $s$ is replaced with, say, $3s$) we may assume that $\alpha^*(G) < s$ and $\alpha(G) < s$.
	If $s \ge n^{1/19}$, then the theorem holds since
	$$
	f(G) \ge \sqrt{\chi(G)} \ge \sqrt{n/\alpha(G)} \ge n^{1/2}s^{-1/2} \ge n^{5/9}s^{-14/9}\;.
	$$
	So assume $s < n^{1/19}$. Let $k = n^{4/9}s^{14/9}$ and observe that $k \ge s^7$. If there are fewer than $(k/se)^k$ acyclic $k$-sets, then by Corollary \ref{coro:prob},
	$f(G) \ge n/k = n^{5/9}s^{-14/9}$ and we are done.
	Otherwise, by Lemma \ref{l:l22}, with $q=n' s^6/k^{1-o(1)}$, $G$ has a $q$-almost-acyclic subgraph $G'$ on
	$n' \ge k^{2-o(1)}/s$ vertices. We claim that $q \ge (n')^{1/3}s^{-2/3}$,
	or equivalently, $(n')^{2/3}/k^{1-o(1)} \ge s^{-20/3}$.
	Indeed,
	$$
	\frac{(n')^{2/3}}{k^{1-o(1)}} \ge \frac{k^{4/3-o(1)}}{s^{2/3}k^{1-o(1)}} \ge \frac{k^{1/3-o(1)}}{s^{2/3}} \ge s^{-20/3}\;.
	$$
	We further claim that $q \le (n')^{1-\Omega(1)}$, or equivalently, $(n')^{-\Omega(1)} \ge s^6/k^{1-o(1)}$.
	Indeed since $k \ge n^{4/9}$, 
	$$
	\frac{s^6}{k^{1-o(1)}} \le \frac{k^{6/7}}{k^{1-o(1)}} = k^{-1/7+o(1)} \le n^{-\Omega(1)} \le 
	(n')^{-\Omega(1)}\;.
	$$
	So, the conditions of Lemma \ref{l:l23} are met, hence $G'$ has an acyclic subgraph $G''$ with
	\begin{align*}
	\chi(G'')  & \ge (n')^{3/4-o(1)}(n' s^6/k^{1-o(1)})^{-1/4}s^{-3/2}\\
	&  = (n')^{1/2-o(1)}k^{1/4-o(1)}s^{-3}\\
	& \ge (k^{2-o(1)}/s)^{1/2-o(1)}k^{1/4-o(1)}s^{-3} \\
	& \ge k^{5/4-o(1)}s^{-7/2}\\
	& = (n^{4/9}s^{14/9})^{5/4-o(1)}s^{-7/2}\\
	& = n^{5/9-o(1)}s^{-14/9}\;.
	\end{align*}
\end{proof}

It remains to prove Lemmas \ref{l:l22} and \ref{l:l23}. Recall that they are respectively analogous to Lemmas 2.2 and 2.3 in \cite{FKS-2021} and their proof flow are similar as well. However, since they each require extensive modifications, we prove them in full detail.
The proof of Lemma \ref{l:l22} which is the key technical lemma is deferred to Section \ref{sec:l22}.
For the proof of Lemma \ref{l:l23} we first require the following simple observation.

\begin{lemma}\label{l:l27}
	Let $G$ be an $n$-vertex orientation with $\alpha(G) \le s$. For a randomly chosen $\pi$, with probability at least
	$1-2^{-\sqrt{n}}$ we have $\alpha(G_\pi) \le 4\sqrt{ns}$.
\end{lemma}
\begin{proof}
	Let $a_k$ be the number of acyclic $k$-sets of $G$.
	By Lemma \ref{l:prob}, for any $k$, the expected number of independent sets of size $k$ in $G_\pi$  is at most
	$$
	a_k \left(\frac{se}{k}\right)^k \le \binom{n}{k}\left(\frac{se}{k}\right)^k
	\le \left(\frac{ne}{k}\right)^k\left(\frac{se}{k}\right)^k = \left(\frac{nse^2}{k^2}\right)^k\;.
	$$
	If $k=4\sqrt{ns}$, then the r.h.s. of the last inequality is at most $2^{-k} \le 2^{-\sqrt{n}}$.
\end{proof}

\begin{proof}[Proof of Lemma \ref{l:l23}]
Let $t=\sqrt{n/q}=n^{\Omega(1)}$.
Consider a permutation $\sigma$ demonstrating the $q$-almost-acyclicity of $G$.
For each $v$, let $B(v)$ be the set of vertices such that $x \in B(v)$ if and only if there is an edge of $G$ 
between 
$v$ and $x$ whose orientation is inconsistent with $\sigma$. Notice that $|B(v)| \le q$.
Partition $\sigma$ into $t$ contiguous blocks $A_1,\ldots,A_t$ of size $n/t$ each.
Now, randomly and independently permute each block to obtain a (semi)random permutation $\pi$.
Since $\alpha(G[A_i]) \le \alpha(G) \le s$, 
by Lemma \ref{l:l27} and the union bound, we may and will assume that $\pi$ is such that for
each $1 \le i \le t$ and each vertex $v$,
\begin{equation}\label{e:1}
	\alpha(G_\pi[A_i])=O(\sqrt{ns/t})\;,
\end{equation}
\begin{equation}\label{e:2}
	\alpha(G_\pi[A_i \cap B(v)])=n^{o(1)}+ O(\sqrt{|A_i \cap B(v)|s})\;.
\end{equation}
Indeed, to see \eqref{e:2}, notice that it trivially holds if $|A_i \cap B(v)| \le \log^4(n) = n^{o(1)}$. Otherwise, by Lemma \ref{l:l27} the failure probability is at most
$2^{-\sqrt{|A_i \cap B(v)|}} \le 2^{-\log^2 n} \ll  1/n^2$ and there are at most $n^2$ events to consider.
Let $Q$ be an independent set in $G_\pi$. We shall bound the cardinality of $Q$.
Select a subset of vertices of $Q$ according to the permutation of $Q$ induced by $\pi$,
as follows. Let $v_1$ be the first vertex of $Q$. Now, assuming we have already selected $v_1,\ldots,v_j$,
let $v_{j+1}$ be the smallest vertex of $Q$ which is not adjacent in $G$ to any vertex $v_1,\ldots,v_j$.
If there is no such vertex, we halt. Denoting the chosen vertices by $v_1,\ldots,v_r$, we have that
they form an independent set in $G$, so $r \le s$, and they are monotone increasing in $\pi$. Let $A_{i_j}$ be the block containing $v_j$. Observe that $i_1 \le i_2 \le \cdots \le i_r$.

We claim that $Q \subseteq \cup_{j=1}^r (A_{i_j} \cup B(v_j))$. Indeed, suppose $x \in Q$. If $x$ is one of the $v_j$'s, the claim is obvious. Otherwise, $x$ must be adjacent in $G$ to one of the vertices $v_1,\ldots,v_r$.
Let $j$ be the smallest index such that $x$ is adjacent in $G$ to $v_j$. Since $x$ is not adjacent to
$v_j$ in $G_\pi$ (as $Q$ is independent in $G_\pi$) and since $x$ is after $v_j$ in $\pi$ (otherwise $x$ would have been chosen, and not $v_j$), we must have that $x$ points to $v_j$ in $G$,
so $x \in A_{i_j} \cup B(v_j)$.

By \eqref{e:1}, we have that $|Q \cap A_{i_j}|=O(\sqrt{ns/t})$, so
$|Q \cap \cup_{j=1}^r A_{i_j}| = O(r\sqrt{ns/t})$. The remaining elements of $Q$ must be in sets of the form
$B(v_j) \cap A_i$ where $A_i$ is any block which contains at least one vertex of $B(v_j)$.
By \eqref{e:2}, the number of elements of $Q$ in each such set is at most
$n^{o(1)}+ O(\sqrt{|A_i \cap B(v_j)|s})$. We have:
\begin{align*}
	|Q| & \le O(r\sqrt{ns/t}) + \sum_{j=1}^r \sum_{i} \left(n^{o(1)}+ O(\sqrt{|A_i \cap B(v_j)|s}\right)\\
	& \le O(r\sqrt{ns/t}) + rtn^{o(1)} + \sum_{j=1}^r O(\sqrt{qts}) \\
	& = O\left(r\sqrt{\sqrt{n}\sqrt{q}s}\right) + rq^{-1/2}n^{1/2+o(1)}\\
	& \le O(s^{3/2}q^{1/4}n^{1/4})+ sq^{-1/2}n^{1/2+o(1)}\\
	& \le s^{3/2}q^{1/4}n^{1/4+o(1)}
\end{align*}
where in the second line we have used that fact that $|B(v_j)| \le q$ and the concavity of the square root function and in the last line we have used our assumption that $q \ge s^{-2/3}n^{1/3}$.
\end{proof}

\section{Proof of Lemma \ref{l:l22}}\label{sec:l22}

As shown in \cite{FKS-2021}, the proof of Lemma \ref{l:l22} follows easily by repeated applications of the following lemma \ref{l:l24}. Before presenting the latter, we require a definition. An orientation is {\em $(q,r)$-almost acyclic} if one can delete from it at most $r$ vertices to make it $q$-almost-acyclic.

\begin{lemma}\label{l:l24}
	Let $n$ be sufficiently large and let $G$ be an $n$-vertex orientation with $\alpha(G) < s$ and $\alpha^*(G) < s$ having at least
	$M$ acyclic $k$-sets. If $q = 6000s^6 n\log^2n/k$, then there is $k'$ satisfying $k-3k/\log n \le k' \le k$ such that $G$ has an induced $(q,6sk/\log n)$-almost-acyclic subgraph $G'$ on at least $e^{-8}kM^{1/k}$ vertices having at least $e^{-6k}M$ acyclic $k'$-sets.
\end{lemma}
\begin{proof}
	Let $t$ be a positive integer parameter to be set later.
	For each $i=0,\ldots,t$ we define vertex sets $W_i,V_i,V'_i$ such that
	\begin{enumerate}
		\item
		$W_i$ is an ordered acyclic $i$-set, $W_i \subset W_{i+1}$, and the ordering of $W_{i+1}$ is consistent with the ordering of $W_i$.
		\item
		$V_0'=V(G)$, $V_i' \supseteq V_i \supset V_{i+1}' \supseteq V_{i+1}$ for $i=0,\ldots,t-1$.
		\item
		$V_i' \cap W_i = \emptyset$.
	\end{enumerate}
	Describing the iterative procedure building these sets requires several notions.
	\begin{itemize}
		\item
		Let $H'_i$ be the $(k-i)$-uniform hypergraph with vertex set $V'_i$ whose edges are the
		$(k-i)$-subsets $S$ for which $G[S \cup W_i]$
		is an acyclic $k$-set for which there is a topological sort $\sigma_S$ consistent with the ordering of $W_i$.
		Let $N'_i$ denote the number of edges of $H'_i$ and observe that $N'_0 \ge M$.
		\item
		Let $H_i$ be an induced subhypergraph of $H'_i$ with minimum degree at least $N'_i/|V'_i|$ (it will always be the case that $N'_i > 0$). Such a subhypergraph can be obtained by repeatedly removing vertices with degree smaller than
		$N'_i/|V'_i|$ until none are left. Let $N_i$ be the number of edges of $H_i$
		and let $V_i$ be its vertex set.
		\item
		Let the ordering of $W_i$ be $\{w_1,\ldots,w_i\}$ so that if $(w_j,w_{j'}) \in E(G)$, then $j < j'$.
		Consider some vertex $v \in V$ and some edge $S$ of $H_i$ containing $v$.
		Recall that $G[S \cup W_i]$ is an acyclic $k$-set for which $\sigma_S$ is a topological sort consistent with the ordering of $W_i$. So if $(w_j,v) \in E(G)$ and $(v, w_{j'}) \in E(G)$, it must be that $j < j'$. Thus, all out-neighbors of $v$ in $W_i$ (if there are any) appear in $W_i$'s ordering after all in-neighbors of $v$ in $W_i$ (if there are any). We say that the set of {\em allowed bins} of $v$ with respect to $W_i$ is $\{j,\ldots,j'\}$
		if $(w_j,v) \in E(G)$, $(v,w_{j'+1}) \in E(G)$ and $v$ is not adjacent in $G$ to any of
		$w_{j+1},\ldots,w_{j'}$. If $v$ has no in-neighbors in $W_i$, then let $j=0$ and if $v$ has no out-neighbors in $W_i$, then let $j'=i$. Notice that for any edge $S$ of $H_i$ containing $v$,
		$\sigma_S$ places $v$ after some $w_\ell$ and before some $w_{\ell+1}$ where $j \le \ell \le j'$ (if $\ell = 0$, then $j=0$ so $v$ is placed before $w_1$ and if $\ell =i$, then $j'=i$ so $v$ is placed after all $w_i$). We say that bin $\ell$ is {\em associated with $(\sigma_S,v)$}. One may think of $\{w_1,\ldots,w_i\}$ as labeling the consecutive internal walls of a set of $i+1$ bins labeled $0,\ldots,i$.
		\item
		Two vertices $x,y$ of $V_i$ are {\em incompatible} if $(x,y) \in E(G)$ and every allowed bin of $x$ with respect to $W_i$ is {\em after} every allowed bin of $y$ with respect to $W_i$
		(in particular, no edge of $H_i$ contains both $x$ and $y$).
	\end{itemize}
	For $0 \le i \le t-1$, we shall go from step $i$ to step $i+1$ by picking some (carefully selected) vertex $w \in V_i$ and some (carefully selected) allowed bin $\ell$ of $w$. We add $v$ to $W_i$ to obtain $W_{i+1}$ with the ordering of $W_{i+1}$ being $\{w_1,\ldots,w_\ell,w,w_{\ell+1},\ldots,w_i\}$.
	We remove from $V_i$ all vertices incompatible with $w$ to obtain $V'_{i+1}$. The edges of $H'_{i+1}$ are obtained by taking all edges $S$ of $H_i$ which contain $w$ and for which $(\sigma_S,w)$ is associated with $\ell$, and keeping $S \setminus w$ as an edge of $H'_{i+1}$.
	
	We now describe the selection process of $w$ and its allowed bin. We do so using two procedures.
	The first procedure, {\em refinement}, is preformed during the first $z := k/\log n$ steps, so at the end of the refinement phase we have $W_z, V'_z, V_z$.
	The second procedure, {\em alignment}, is preformed during the final $t-z$ steps, so at the end of the alignment phase we obtain $W_t,V'_t,V_t$.
	
	\begin{itemize}
		\item
		{\bf Refinement.} This is the more technical procedure and requires some setup.
		Each vertex is designated either {\em white} or {\em black}, where in the beginning (at step $i=0$)
		all vertices are white. Once a vertex becomes black it remains so.
		We shall keep the number of black vertices after step $i$ at most $2s i$.
		The property required of all white vertices is that their number of
		allowed bins is at most $2s$.
		For a bin $0 \le j \le i$, let $X_{i,j}$ be the set of white vertices that have bin $j$ allowed. Let $p_i = \arg\max_{0 \le j \le i}|X_{i,j}|$ be a bin most allowed for white vertices.
		
		For $v \in V_i$, let $p_{i,v}$ be an allowed bin for $v$ associated with maximum number of $(\sigma_S,v)$ where $S$ is an edge of $H_i$ containing $v$.
		Consider some $v \in X_{i,p_i}$ and notice that since $v$ is white, it has at most $2s$ allowed bins so $|p_i-p_{i,v}| \le 2s-1$.
		Hence, there is a set $Y \subseteq X_{i,p_i}$ of size at least $|X_{i,p_i}|/4s$ such that all $v \in Y$
		share the same $p_{i,v}$, call it $p^*_i$. We thus have that
		$$
		|Y| \ge \frac{|X_{i,p_i}|}{4s}\;.
		$$
		
		By Lemma \ref{l:deg}, $G[Y]$ has a vertex $w$ whose in-degree and out-degree in $G[Y]$ is at least $|Y|/4s-1/2$. This is our chosen vertex $w$ and $p^*_i$ is the chosen allowed bin.
		
		Adding $w$ to $W_i$ so as to obtain $W_{i+1}$ corresponds to splitting bin $p^*_i$ into two bins where the new internal wall is labeled by $w$.
		However, it might be that some white vertices which have $2s$ allowed bins prior to the split,
		now have $2s+1$ allowed bins and must therefore be labeled black. We claim that there are fewer than $2s$ such problematic vertices. Indeed, each problematic vertex either has $s$
		allowed bins to the left of the new internal wall $w$, or else has $s$ allowed bins to the right of the new internal wall $w$. So, if there were $2s$ problematic vertices, without loss of generality, at least $s$ of them have $s$ allowed bins to the left of the new internal wall $w$. But this means that in $G$ they are not adjacent to any of
		$w,w_{p^*_i-1},w_{p^*_i-2},\ldots,w_{p^*_i-s+1}$ contradicting the assumption that
		$\alpha^*(G) < s$.
		Thus, we have fewer than $2s$ problematic vertices. Designating them black
		increases the number of black vertices by at most $2s$ so there are now at most $2si+2s=2s(i+1)$  black vertices at step $i+1$.

		\item
		{\bf Alignment.} Let $\varepsilon = \log^2 n/k$. Take a white vertex $w$ incompatible with at least $\varepsilon|V_i|$ other vertices of $V_i$ (observe that alignment may not be possible)
		and let the allowed bin be $p_{i,w}$ where recall that $p_{i,w}$ is a bin with maximum number of associated $(\sigma_S,w)$. Split bin $p_{i,w}$ precisely as in refinement (namely, $p_{i,w}$ plays the role of $p^*_i$ in refinement).
		Make black all vertices that had $2s$ allowed bins and now have $2s+1$ allowed bins.
		We have already shown that the number of new black vertices is at most $2s$, so we keep the required upper bound of $2si+2s = 2s(i+1)$ on the amount of black vertices at step $i+1$.
	\end{itemize}
	
	We halt when it is no longer possible to do an alignment; this defines $t$.
	Observe that for an alignment step we have $|V'_{i+1}| \le (1-\varepsilon)|V_i|
	\le (1-\varepsilon)|V'_i|$, so we cannot perform more than $\log_{1/(1-\varepsilon)} n < (2/\varepsilon) \log n = 2k/\log n -1 $ alignment steps, hence $t \le z + 2k/\log n -1 = 3k/\log n -1$.
	
	\begin{claim}\label{cl:1}
		$N_t \ge e^{-6k}M$ and $|V_t| \ge e^{-8}kM^{1/k}$.
	\end{claim}
	\begin{proof}
		Observe that the edges of $H'_{i+1}$ are obtained by taking all edges $S$ containing a white vertex
		$w$ such that $(\sigma_S,w)$ is associated with the most popular allowed bin of $w$ and deleting $w$ from each of them. Let $d_i(w)$ be the degree of $w$ in $H_i$ and recall that 
		$d_i(w) \ge N'_i/|V'_i|$. Hence, $N'_{i+1} \ge N'_i/r|V'_i|$ where $r$ is the number of allowed bins for $w$ with respect to $W_i$.
		As $w$ is white, the number of allowed bins is at most $2s$. Whence, $N'_{i+1} \ge N'_i/2s|V'_i| \ge N'_i/n^2$, so
		$$
		N_t \ge  \frac{N'_t}{n^2} \ge \frac{N'_0}{(n^2)^{t+1}} \ge \frac{M}{(n^2)^{t+1}} \ge \frac{M}{(n^2)^{3k/\log n}} \ge e^{-6k}M\,.
		$$
		Since $H_t$ is a $(k-t)$-uniform hypergraph, its number of edges $N_t$ trivially satisfies 
		$$
		N_t \le \binom{|V_t|}{k-t} \le \left(\frac{e|V_t|}{k-t} \right)^{k-t}
		$$
		implying that
		$$
		|V_t| \ge \left(\frac{k-t}{e} \right)N_t^{1/(k-t)} \ge \left(\frac{k-t}{e} \right) (e^{-6k}M)^{1/(k-t)}
		\ge e^{-8}kM^{1/k}\;.
		$$
	\end{proof}
	\begin{claim}\label{cl:2}
		After step $z$ (i.e., at the end of the refinements), we have $X_{z,j} \le 2560s^5 n \log n /z$ for all $0 \le j \le z$. Consequently, $|X_{t,j}| \le 2560s^5 n \log n/z$ for all $0 \le j \le t$.
	\end{claim}
	\begin{proof}
		It suffices to prove that $X_{z,j} \le 2560s^5n \log n/z$ for all $0 \le j \le z$ as the amount of white vertices in a bin can never grow over time.
		
		Referring back to the refinement procedure at step $i$, recall that $|X_{i,p_i}|$ is maximum
		among all other $|X_{i,j}|$. Further, we designate in it a subset $Y$ of size
		$|Y| \ge |X_{i,p_i}|/4s$ such that all vertices $v \in Y$ 
		share the same $p_{i,v}$ which we have denoted by $p^*_i$.
		Also recall that $w$ is a vertex of $Y$ such that in $G[Y]$, the in-degree and out-degree of $w$ is at least $|Y|/4s-1/2 \ge |Y|/5s$ (since $|Y| \ge 10s$; indeed, we may assume $|X_{i,p_i}| \ge 2560s^5 n \log n/z > 2560s^5$ otherwise there is nothing to prove, so $|Y| \ge 640s^4 \ge 10s$).
		We then split bin $p^*_i$ in two. Immediately prior to the split, bin $p^*_i$ has at least
		$|Y|$ {\em balls} (for a white vertex $v$ and an allowed bin for it, view a ball labeled $v$ in that bin). Immediately after the split, the left split bin has a ball for all in-neighbors of $Y$ in $X_{i,p^*_i}$, in particular, the in-neighbors of $Y$ in $G[Y]$, so has at least $|Y|/5s \ge |X_{i,p_i}|/20s^2$ balls.
		Similarly, the right split bin has a ball for all out-neighbors of $Y$ in $X_{i,p^*_i}$, in
		particular, then out-neighbors of $Y$ in $G[Y]$, so has at least $|Y|/5s \ge |X_{i,p_i}|/20s^2$ balls.
		
		We can view the evolution of the white vertices in the bins using a full binary tree $T$.
		Each node contains the set of white vertices (balls) at its creation time, so the root (which corresponds to the single bin at step $i=0$) contains all of $V$. Each internal node (i.e., non-leaf) is labeled with a pair $(i,j)$ meaning that it represents the unique bin $j$ at step $i$ that was split (so the root is labeled $(0,0)$). Each internal node thus has two children. The leaves represent the final $z+1$ bins
		obtained after the final refinement. Again, we emphasize that the subset of vertices (balls) of an internal node labeled $(i,j)$ is precisely the set $X_{i,j}$ at the time of the split.
		As $T$ is full binary with $z+1$ leaves, it has at least $(z+1)/2$ leaf parents.
		Consider the leaf parent labeled $(z-1,j)$ (namely, the final refinement step splits some bin $j$ into two bins). Suppose that the $|X_{z-1,j}|=r$. So we have that
		$r \le |X_{z-1,p_{z-1}}| \le 4sr$
		and recall that $|X_{z-1,p_{z-1}}|$ is the largest at step $z-1$. Since
		$|X_{i,p_i}| \le |X_{i-1,p_{i-1}}|$ for all $1 \le i \le z-1$ (the number of balls in bins can only decrease so the maximum cannot increase over time) and since the number of balls in a parent is never smaller than the number of balls in a child, this means that the number of balls in all internal vertices is at least $r/4s$.
		How many elements are in $X_{z,j'}$ for any $0 \le j' \le z$? If there were more than $4sr$, then its parent would also have more than $4sr$ balls, but then the split bin at stage $z-1$ would have more than $r$ balls, contradicting our assumption that $|X_{z-1,j}|=r$. So, we have that $|X_{z,j'}| \le 4sr$.
		
		Assume, to the contrary, that $|X_{z,j}| > 2560s^5n \log n/z$ for some $j$. Then
		$r \ge 640 s^4 n \log n/z$. But since there are at least $z/2$ internal vertices, the number of balls
		is all of them is at least $80s^3n \log n$. As there are no more than $n$ vertices, there is some vertex with a ball representing it in at least $80s^3 \log n$ internal nodes.
		But now let us consider the height of $T$ which, we claim, is not larger than
		$40 s^2 \log n$.
		Consider some internal node labeled $(i,j)$. Since
		$|X_{i,j}| \le |X_{i,p_i}|$ and since we have already seen that each child loses at least
		$|X_{i,p_i}|/20s^2$ balls (that are only present in its sibling), the ratio of the number of balls between a child and a parent is at most $(1-1/20s^2)$.
		This means that the height of the tree is not larger than
		$\log_{(1/(1-1/20s^2))} n < 40s^2 \log n$.
		Now, consider all the $80s^3 \log n$ nodes having a ball representing the same vertex $v$.
		Construct a directed acyclic graph $D$ whose vertices are these nodes and there is a directed edge from a node $x$ to a node $y$ if $y$ is an ancestor of $v$. So the maximum out-degree
		of $D$ is smaller $40s^2 \log n$. But this means that $D$ has an independent set of size larger than $2s$. However, this is impossible since consider the node of this independent set
		having label $(i,j)$ where $i$ is maximum among all the nodes of the independent set.
		This precisely means that at step $i$, $v$ could not have been white as it has more than $2s$ allowed bins, a contradiction.
	\end{proof}
	
	Let $G'=G[V_t]$. By Claim \ref{cl:1}, $G'$ has at least $e^{-8}kM^{1/k}$ vertices and at least
	$e^{-6k}M$ acyclic $k'$-sets for some $k' \ge k-3k/\log n$. It remains to show that
	for $q = 6000s^6 n\log^2n/k$,  $G'$ is $(q,6sk/\log n)$-almost-acyclic.
	Recall that the number of black vertices of $V_t$ is at most $2st < 6sk/\log n$.
	So, it remains to show that if we remove from $G'$ all the black vertices, we remain with
	a subgraph $G''$ that is $q$-almost acyclic (note: we are not claiming that the number of acyclic $k'$-sets in $G''$ is large). For $v \in V_t$, let $b_v$ be the first allowed bin for $v$, and for a bin $j$, let $B_j$
	be the set of all white vertices $v \in V_t$ with $b_v=j$.
	Consider an arbitrary permutation of the white vertices of $V_t$ such that for $j < j'$,
	all the elements of $B_j$ appear before all the elements of $B_{j'}$.
	For any such white $v$, there are at most $\varepsilon |V_t| < \epsilon n$ vertices incompatible with $v$.
	Furthermore if $b_v=j$, the last allowed bin for $v$ is at most $j+2s-1$ so all white vertices
	$v'$ with $b_v' \ge j+2s$ are necessarily compatible with $v$.
	But how many vertices could we have in any $B_j$? By Claim \ref{cl:2} there are at most
	$2560s^5 n \log n/z$ such vertices. Hence there could additionally be at most $5120s^6 n \log n/z$
	edges incident with $v$ which are oriented against our chosen permutation. So, $G''$ is $q'$-almost-acyclic already for
	$$
	q' = \varepsilon n + 5120s^6 n \log n/z \le \varepsilon n + 5120s^6 n \log^2 n/k < 5121s^6 n  \log^2 n/k < q\;.
	$$
\end{proof}

\begin{proof}[Proof of Lemma \ref{l:l22}]
	Let $t = \sqrt{\log n}$. Construct a sequence of subgraphs $G = G_0 \supseteq G_1 \supseteq \cdots \supseteq G_t$  iteratively applying Lemma \ref{l:l24}.
	For each $0 \le i \le t$ denote by $n_i$ the number of vertices in $G_i$ and let $k_i$ be the
	value of $k'$ in Lemma \ref{l:l24} for which we obtain a bound on the number of acyclic $k_i$-sets in $G_i$.
	Notice that $k_0 = k$ and for $1 \le i \le t$ it holds that  $k_i \ge k_{i-1}(1-3/\log n_{i-1})$. 
	Let $q_i = 6000s^6 n_{i-1}\log^2n_{i-1}/k_{i-1}$. Then, each $G_i$ is
	$(q_i,6sk_{i-1}/\log n_{i-1})$-almost-acyclic, has at least $M_i := e^{-6(k_0+\cdots+k_{i-1})}(k/se)^k \ge e^{-6ik}(k/se)^k$
	acyclic $k_i$-sets, and for $i > 0$, has
	$$
	n_i \ge e^{-8}k_{i-1}(M_{i-1})^{1/k_{i-1}} \ge e^{-8}k_{i-1}(e^{-6ik}(k/se)^k)^{1/k} = \Omega(e^{-6i}k_{i-1}k/s)
	$$
	vertices. Recalling that $k_i \ge k_{i-1}(1-3/\log n_{i-1})$ and that $k = n^{\Omega(1)}$, one immediately obtains by induction that $n_i \ge k^{2-o(1)}/s$ and $k_i \ge ke^{-O(i/\log n)} = (1-o(1))k$ for every $i \le t = \sqrt{\log n}$. Hence, $q_i = s^6n_{i-1}k^{o(1)-1}$.
	As trivially there is some $i$ such that $\log n_i-\log n_{i-1} = \log(n_i/n_{i-1}) \le (1/t)\log n \le \sqrt{\log n}$, we have that $n_i/n_{i-1} = k^{o(1)}$. Then $G'=G_i$ satisfies the
	properties stated in the lemma once we remove from it at most $6sk_{i-1}/\log n_{i-1} < sk$ vertices. But notice that as $G_i$ has at least $k^{2-o(1)}/s$ vertices and since $s \le k^{1/7}$, the fraction of vertices removed is less than half. So, after removal we obtain a graph which still
	has $n'=k^{2-o(1)}/s$ vertices and is at least $s^6n'k^{o(1)-1}$-almost-acyclic.
\end{proof}

\section*{Acknowledgment}
I thank referees $X$ and $Y$ for very useful comments and very careful reading of the manuscript.
Special thanks to referee $Y$ who spotted a gap in the proof of Lemma \ref{l:l24} in an earlier version of this paper.



\end{document}